\documentclass[12pt]{amsart}
\usepackage{amssymb,amsmath,amsthm,comment}
\usepackage{pdfpages,graphicx} 
\usepackage[section]{placeins} % figures should be displayed in same section
\usepackage{wrapfig}
\usepackage{float}
\usepackage{lineno}
\usepackage{youngtab}

\usepackage{setspace}
\newcommand{\shrinkmargins}[1]{
  \addtolength{\textheight}{#1\topmargin}
  \addtolength{\textheight}{#1\topmargin}
  \addtolength{\textwidth}{#1\oddsidemargin}
  \addtolength{\textwidth}{#1\evensidemargin}
  \addtolength{\topmargin}{-#1\topmargin}
  \addtolength{\oddsidemargin}{-#1\oddsidemargin}
  \addtolength{\evensidemargin}{-#1\evensidemargin}
  }
\shrinkmargins{0.5}

\newtheorem{theorem}{Theorem}
\newtheorem{lemma}[theorem]{Lemma}

\newtheorem*{theorem*}{Theorem}

\newtheorem{proposition}[theorem]{Proposition}
{Claim}

\theoremstyle{definition}
\newtheorem*{definition}{Definition}

\theoremstyle{remark}

\newtheorem*{remarks}{{\bf Remarks}}

\numberwithin{theorem}{section} \numberwithin{equation}{section}

\usepackage{caption}
\usepackage{multirow}

\def\func#1{\mathop{\rm #1}}%
\begin{document}
\title[Tur\'an Inequalities]{Tur\'an Inequalities for Infinite Product Generating Functions}
\author{Bernhard Heim }
\address{Lehrstuhl A f\"{u}r Mathematik, RWTH Aachen University, 52056 Aachen, Germany}
\email{bernhard.heim@rwth-aachen.de}
\author{Markus Neuhauser}
\address{Kutaisi International University, 5/7, Youth Avenue,  Kutaisi, 4600 Georgia}
\address{Lehrstuhl A f\"{u}r Mathematik, RWTH Aachen University, 52056 Aachen, Germany}
\email{markus.neuhauser@kiu.edu.ge}
\subjclass[2010] {Primary 05A17, 11P82; Secondary 05A20}
\keywords{Generating functions, Log-concavity, Tur\'an inequalities, Partition numbers.}
%%\linenumbers
\begin{abstract}
In the $1970$s, Nicolas proved that the partition function $p(n)$ is log-concave for
$ n > 25$. In \cite{HNT21}, a precise conjecture on the log-concavity for 
the plane partition function $\func{pp}(n)$ for $n >11$ was 
stated. This was recently proven by 
Ono, Pujahari, and Rolen. 
In this paper, we provide a general picture. 
We associate to double sequences 
$\{g_d(n)\}_{d,n}$ with $g_d(1)=1$ and
$$0 \leq g_{d}\left( n\right) -n^{d}\leq
g_{1}\left( n\right)
\left( n-1\right) ^{d-1}$$
polynomials $\{P_n^{g_d}(x)\}_{d,n}$ given by
\begin{equation*}
\sum_{n=0}^{\infty} P_n^{g_d}(x) \, q^n := 
\func{exp}\left( x \sum_{n=1}^{\infty} g_d(n) \frac{q^n}{n} \right)
=\prod_{n=1}^{\infty} 
\left( 1 - q^n \right)^{-x f_d(n)}.
\end{equation*}
We recover $ p(n)= P_n^{\sigma_1}(1)$ 
and $\func{pp}\left( n\right) = P_n^{\sigma_2}(1)$, 
where $\sigma_d (n):= \sum_{\ell 
\mid n} \ell^d$ and $f_d(n)= n^{d-1}$.
Let $n \geq 6$.
Then the sequence $\{P_n^{\sigma_d}(1)\}_d$ is log-concave for
almost all $d$ if and only if $n$ is divisible by $3$.
Let $\func{id}(n)=n$. Then $P_n^{\func{id}}(x) = \frac{x}{n} L_{n-1}^{(1)}(-x)$, where
$L_{n}^{\left( \alpha \right) }\left( x\right) $
denotes the $\alpha$-associated Laguerre polynomial.
In this paper, we invest in Tur\'an inequalities 
\begin{equation*}
\Delta_{n}^{g_d}(x) := \left( P_n^{g_d}(x) \right)^2 -  
P_{n-1}^{g_d}(x)  \,  P_{n+1}^{g_d}(x)  \geq 0.
\end{equation*}
Let $n \geq 6$ and $0 
\leq x < 2 - \frac{12}{n+4}$.
Then
$n$ is divisible by $3$ if and only if $\Delta_{n}^{g_d}(x) \geq 0$ for
almost all $d$. Let $n \geq 6$ and $n \not\equiv 2 \pmod{3}$. 
Then the condition on
$x$ can be reduced to $x 
\geq 0$.
We determine explicit bounds. As an analogue
to Nicolas' result,
we have for $g_1= \func{id}$ that
$\Delta_{n}^{\func{id}}(x)  \geq 0$ for all
$x \geq 0 $ and all $n$. 
\end{abstract}
%%%%%%%%%%%%%%%%%%%%%%%%%%%%%%%%%%%%%%%%%%%%%%%%%%%%%%%%%%%%%%%%%%%%%%%%%%%%%%%%%%%%%%%%%%%%%%%%%%%%%%%%%%
%%%%%%%%%%%%%%%%%%%%%%%%%%%%%%%%%%%%%%%%%%%%%%%%%%%%%%%%%%%%%%%%%%%%%%%%%%%%%%%%%%%%%%%%% Section 1
%%%%%%%%%%%%%%%%%%%%%%%%%%%%%%%%%%%%%%%%%%%%%%%%%%%%%%%%%%%%%%%%%%%%%%%%%%%%%%%%%%%%%%%%%%%%%%%%%%%%%%%%%%
%%%%%%%%%%%%%%%%%%%%%%%%%%%%%%%%%%%%%%%%%%%%%%%%%%%%%%%%%%%%%%%%%%%%%%%%%%%%%%%%%%%%%%%%%%%%%%%%%%%%%%%%%%
\maketitle
\newpage
\section{Introduction and main results}
In this paper, we study Tur\'an inequalities $p_n(x)^2
-p_{n-1}(x) \, p_{n+1}(x) \geq 0$
for families of polynomials $\{p_n(x)\}_n$ attached to arithmetic functions.

Our work is motivated by a recent result by Ono, Pujahari,
and Rolen \cite{OPR22} and \cite{HN22}. 
Ono, Pujahari, and Rolen proved
the log-concavity conjecture (\cite{HNT21}, Conjecture 1)
for plane partitions
$\func{pp}\left( n\right) $ for $n >11$.
Twenty-four years ago, Nicolas \cite{Ni78} had already proved the 
log-concavity property for the partition numbers $p(n)$ for $n >25$. This result was reproved by DeSalvo and Pak \cite{DP15}.
For an introduction %in
to partition numbers and plane partition numbers, we refer to Andrews'
book \cite{An98}. Further, to study
the concept of log-concavity and related topics, Brenti \cite{Br89} and Stanley \cite{St89,St99} are suitable references.

This paper is also a significant generalization of our previous result in \cite{HN22}.
Let $a_d(n):= \frac{1}{n} \sum_{k=1}^n \left( \sum_{\ell \vert k} \ell^d \right) \, a_d(n-k)$, with $a_d(0)=1$. Let $n \geq 6$ be fixed.
Then the sequence $\left\{ a_d\left( {n}\right) \right\}_{d}$ is log-concave for almost all $\in \mathbb{N}$ if and only
if $n \equiv 0 \pmod{3}$.
Note that $p(n)= a_1(n)$ and $\func{pp}(n)= a_2(n)$.
The quantities $p(n)$ and $\func{pp}\left( n\right) $ are 
induced by certain arithmetic functions. %, which    
This leads to the following generalization. 

\begin{definition} Let $\mathbb{D}$ be the set of all double sequences
$\{g_d(n)\}_{d,n \geq 1}$ with normalization $g_d(1)=1$, such that
$\sum_{n=1}^{\infty} g_d(n) \, q^{n-1}$ is regular at $q=0$ with radius of convergence $R$, and
\begin{equation*}\label{condition}
0 \leq g_{d}\left( n\right) -n^{d} \leq 
g_{1} \left( n\right)
\,  \left( n-1\right) ^{d-1}
\end{equation*}
for all $d$ and $n$. 
\end{definition}
We invest in sequences of polynomials $\{P_n^{g_d}(x)\}_n$, defined by the
recurrence relation:
\begin{equation*}
P_n^{g_d}(x) := \frac{x}{n} \, \sum _{k=1}^{n} g_d(k) \, P_{n-k}^{g_d}(x), \qquad \text{with } P_0^{g_d}(x):=1.
\end{equation*}
We have the generating series
\begin{equation*}
\sum_{n=0}^{\infty} P_n^{g_d}(x) \, q^n =
\func{exp}\left( x \sum_{n=1}^{\infty} g_d(n) \frac{q^n}{n} \right) =
 \prod_{n=1}^{\infty} 
\left( 1 - q^n \right)^{-x \,  f_d(n)},
\end{equation*}
where $ n \, f_d(n) = \sum_{\ell \mid n } \mu(\ell) \, g_d(n/ \ell)$ with $\mu$ the Moebius function.
Examples for $g_d(n)$ are $\{\sigma_d(n)\}$ and $\{\psi_d(n)\}$,
where $\sigma_d(n)= \sum_{\ell \mid
n} \ell^d$  and $\psi_d(n)= n^d$. 
Tur\'an's inequality of $\{P_n^{g_d}(x)\}$ at $n$ for a subset of $\mathbb{R}$ is defined by
\begin{equation*}%\label{Turan}
\Delta_{n}^{g_d}(x) := \left( P_n^{g_d}(x) \right)^2 -  
P_{n-1}^{g_d}(x)  \,  P_{n+1}^{g_d}(x)  \geq 0.
\end{equation*}
Let $x_0$ be fixed, we call
$P_n^{g_d}(x_0)$ log-concave at $n$ if $\Delta_{n}^{g_d}(x_0)\geq 0$. 

We note that the partition function and the plane partition function satisfy
$p(n)= P_{n}^{\sigma_1}(1)$ and
$\func{pp}\left( n\right) = P_n^{\sigma_2}(1)$. Let $E^{g_d}$ be the set of all $n \in \mathbb{N}$ with
$\Delta_{n}^{g_d}(1)<0$, denoted as strictly log-convex.

Nicolas \cite{Ni78} proved that the partition function $p(n)$ is log-concave for almost all $n$. The set of exceptions is given by
$E^{\sigma _{1}}= \left\{ %n=
2k +1 \, : \, 0 \leq k \leq 12\right\} $.
Ono, Pujahari, and Rolen \cite{OPR22} proved that $E^{\sigma_2}= \{1,3,5,7,9,11\}$.
Numerical investigations \cite{HN22} for $n \leq 10^5$ indicate that $E^{\sigma_3}= \{1,3,5,7\}$, 
$E^{\sigma_4}= E^{\sigma_5}=\{1,5\}$.
Surprisingly, 
$E^{\sigma_{20}}$ has at least $10$ elements. We believe that the general and {\it clean} patterns 
associated to double sequences in $
\mathbb{D}$ are displayed by $g_d(n)= \psi_d(n)$ (see Table \ref{clean}).
\begin{table}[H]
\[
\begin{array}{r|cccccccccccccccccccc}
\hline
n\backslash d&{1}&{2}&3&4&5&6&7&8&9&10&11&12&13&14&15&16&17&18\\ \hline \hline
1&\bullet &\bullet &\bullet &\bullet &\bullet &\bullet &\bullet &\bullet &\bullet &\bullet &\bullet &\bullet &\bullet &\bullet &\bullet &\bullet &\bullet &\bullet  \\
2&&&&&&&&&&&&&&&&&&\\
3&&&&&&&&&&&&&&&&&&\\
4&&&&&\bullet &\bullet &\bullet &\bullet &\bullet &\bullet &\bullet &\bullet &\bullet &\bullet &\bullet &\bullet &\bullet &\bullet \\
5&&&&&&&&&&&&&&&&&&\\
6&&&&&&&&&&&&&&&&&&\\
7&&&&&&&&&&&\bullet &\bullet &\bullet &\bullet &\bullet &\bullet &\bullet &\bullet  \\
8&&&&&&&&&\bullet &\bullet &\bullet &\bullet &\bullet &\bullet &\bullet &\bullet &\bullet &\bullet \\
9&&&&&&&&&&&&&&&&&&\\
10&&&&&&&&&&&&&&&&\bullet &\bullet &\bullet \\
11&&&&&&&&&&&&\bullet &\bullet &\bullet &\bullet &\bullet &\bullet &\bullet  \\
12&&&&&&&&&&&&&&&&&&\\
13&&&&&&&&&&&&&&&&&&\\
14&&&&&&&&&&&&&&&\bullet &\bullet &\bullet &\bullet  \\
\hline
\end{array}
\]
\caption{\label{clean}
Exceptions for $g_d\left( n\right) =n^{d}$, $1\leq d\leq 18
$ and $1\leq n\leq 14
$
}
\end{table}
In our main result we capture the impact of the residue of $n$ divided by $3$ and
the range of the argument of the $\Delta_n^{g_d}(x)$.
%%%%%%%%%%%%%%%%%%%%%%%%%%%%%%%%%%%%%%%%%%%%%%%%%%%%%%%%%%%%%%%%%%%%%%%%%%%%%%%%%%%%%%%
%%%%%%%%%%%%%%%%%%%%%%%%%%%%%%%%%%%%%%%%%%%%%%%%%%%%%%%%%%%%%%%%%%%%%%%%%%%%%%%%%%%%%%%
%%%%%%%%%%%%%%%%%%%%%%%%%%%%%%%%%%%%%%%    Main Theorem  %%%%%%%%%%%%%%%%%%%%%%%%%%%%%%
%%%%%%%%%%%%%%%%%%%%%%%%%%%%%%%%%%%%%%%%%%%%%%%%%%%%%%%%%%%%%%%%%%%%%%%%%%%%%%%%%%%%%%%
\begin{theorem} \label{th:1}
Let $\{g_d(n)\}$ be a double sequence in $\mathbb{D}$. Let $n \geq 6$. Moreover let
\begin{equation}\label{Turan}
\Delta_{n}^{g_d}(x)= \left( P_n^{g_d}(x) \right)^2 -  P_{n-1}^{g_d}(x) \, P_{n+1}^{g_d}(x) \geq 0
\end{equation}
be the Tur\'an inequality. 
\begin{itemize}
\item[a)] Let $0 \leq x < 2- \frac{12}{n+4}$. Then (\ref{Turan}) %i
holds true for almost all $d$ if and only
if $n$ is divisible by $3$.
\item[b)] Let $n \not\equiv 2 \pmod{3}$ and $x %>
\geq 0$. Then (\ref{Turan}) %i
holds true for almost all $d$ if and
only if $n$ is divisible by $3$.
\end{itemize}
\end{theorem}

%Let
The case $g_{d}\left( n\right) = \sigma_{d}%1
\left( n\right) $ and $x=1$ leads
to the results obtained in \cite{HN22}, Theorem~1.2 and Theorem 1.3.
An explicit analysis of the bounds obtained in the proof of Theorem \ref{th:1} leads to the following:
%%%
%%%
%%%%%%%%%%%%%%%%%%%%%%%%%%%%%%%%%%%%%%%%%%%%%%%%%%%%%  Theorem
\begin{theorem}
Let $\left\{ g_{d}
\left( n\right)
\right\} $ be a double sequence in $
\mathbb{D}$. Let $n \geq 3$ and $n \neq 5$.
Let~$R$ be the radius of convergence of $\sum_{n=1}^{\infty} g_1(n) \, \frac{q^n}{n}$.
For each $x$, let $r(x)$ be chosen with $ 0 < r(x) < R$ and 
$P_n^{g_{1}%d
}(x) \leq r(x)^{-n}$ for all $n$.
Then we have the following properties.
\begin{itemize}
\item[(i)]
Let $n \equiv 0 \pmod{3}$ and $x >0$. Then $\Delta_n^{g_d}(x) \geq 0$ for $d \geq d_{0}\left( n,x\right) $, where
\begin{equation*}
d_{0}\left( n,x\right) =1+\frac{2n}{3\ln \left( 9/8\right) }
\left( \ln \left( n/3\right) -\ln \left( x\right) -3\ln \left( r(x) \right) \right) .
\end{equation*}
\item[(ii)]
Let $n \equiv 1 \pmod{3}$ and $x >0$. Then $\Delta_n^{g_d}(x) < 0$ for $d \geq d_{0%^{\prime
}\left( n,x\right) $, where
\begin{equation*}
d_{0}%^{\prime }
\left( n,x\right) =1+\frac{2n}{3\ln \left( 9/8\right) }
\left( \ln \left( \frac{n-1}{3}\right) %-3\ln \left( r(x) \right)
-\ln \left( x\right) -3\ln \left( r(x) \right) \right) .
\end{equation*}
\item[(iii)]
Let $n \equiv 2 \pmod{3}$ and $0<x<2-\frac{12}{n+4}$. Then
$\Delta_n^{g_d}(x) < 0$ for $d \geq d_{0%^{\prime \prime
}\left( n,x\right) $, where
\begin{eqnarray*}
d_{0}%^{\prime \prime }
\left( n,x\right) & = &
1+\frac{1}{
\ln \left( 9/8\right) }\left( -\ln \left( \frac{n-2}{3n+3}
\left( \frac{1}{x}+\frac{1}{2}\right) -\frac{1}{3}\right)  \right.\\ 
& & %\phantom{1+\frac{1}{\ln \left( 9/8\right) }\left( \right. }
\left. + \, \frac{n-2}{3}\ln \left( \frac{n-2}{3}\right) -
\frac{n+1}{3}\ln \left( x\right) -n\ln \left( r(x) \right) \right) .
\end{eqnarray*}
\end{itemize}
\end{theorem}
%%%%%%%%%%%%%%%%%%%%%%%%%%%%%%%%%%%%%%%%%%%%%%%%%%%%%%%%%%%%%%%%%%%%%  Remarks
\begin{remarks}\ \\
a) The positive real number $r(x)$ always exists due to Cauchy--Hadamard's theorem.\\
b) Let $g_d(n)= n^d$. Then $\Delta_4^{g_5}(x)$ has sign changes for positive real $x$, 
since there are two positive, real zeros $\alpha_1 < \alpha_2$.
\end{remarks}
%%%%%%%%%%
%%%%%%%%%%
%%%%%%%%%%
\section{Records}
Let $g_d(n)= \sigma_d(n)$. Then for $d=1$ and $d=2$ complete results for  
the log-concavity $\Delta_n^{g_d}(1) \geq 0$, including the explicit 
$E^{\sigma_d}$, are provided by Nicolas \cite{Ni78} and
Ono--Pujahari--Rolen \cite{OPR22}.
For $n\leq 10^{5}$ and $d \leq 8$ further results have been obtained by
Heim--Neuhauser \cite{HN22}. 

Let $g_d(n)= \psi_d(n)$. In Table \ref{table:d} we have displayed the results for
$1 \leq d \leq 9$. In this paper, we prove the analogue %y
to Nicolas'
result and give some numerical
evidence for the case $d=2$, which is for $\sigma _{2}\left( n
\right) $ the log-concavity for plane partitions.
\newpage
\begin{table}[H]
\begin{tabular}{r|l|l|l}
\hline
$ d$ & log-concave &  strictly log-convex & proof / verification\\ \hline \hline
$1$ & $n>1$      &  $n=1$   & Heim--Neuhauser 2022             \\
$2$&          $n >1$ &  $n=1$    &     H--N 2022   ($n \leq  10^{4}$)  \\
$3$&      $n >1$ & $n=1$ & \\
%%time = 4h, 5min, 44,594 ms
$4$& $n > 1$ & $n=1$ & 
\phantom{xxxxxx}\vdots
%%%%%%%%%%%%2022 H--N ($n \leq 10^{5
%%%%%%%%%%%%%%%%%%%%%%%%%%%%%%%%%%}$)
%time = 7h, 32,270 ms
\\
$5$&$n >4, \{2,3\}$&           $\{1,4\}$                &  
\phantom{xxxxxx}\vdots
%%2022 H--N ($n \leq 10^{5}$)%time = 11h, 16min, 22,853 ms
\\
$6$&$n >4, \{2,3\}$&$\left\{ 1,4\right\} $& 
\phantom{xxxxxx}\vdots 
%%2022 H--N ($n \leq 10^{3}$)%time = 16h, 18min, 25,314 ms
\\
$7$&$n >4, \{2,3\}$&$\left\{ 1,4\right\} $&  
\phantom{xxxxxx}\vdots 
%% 2022 H--N ($n \leq 10^{3}$)     
\\
$8$&$n >4, \{2,3\}$&$\left\{ 1,4\right\} $&        
\phantom{xxxxxx}\vdots
%%2022 H--N ($n \leq 10^{3}$)
\\
$9$&$n >8 , \{2,3,5,6,7 \}$&$\{1,4,8 \}$&        %2022
H--N 2022 ($n \leq 10^{4}$)\\
\hline
\end{tabular}
\caption{\label{table:d}  Properties of $\{ P_n^{\psi_d}(1)\}$          }
\end{table}
%%%%%%%%%%%%%%%%%%%%%%%%%%%%%%%%%%%%%%%%%%%%%%%%%%%%%%%%%%%%%%%%%%
%%%%%%%%%%%%%%%%%%%%%%%%%%%%%%%%%%%%%%%%%%%%%%%%%%%%%%%%%%%%%%%%%%
More generally, let $\{g_d(n)\}$ be a double sequence in $
\mathbb{D}$.
Then $\Delta_1^{g_d}(1)$ is always negative, since 
\begin{equation*}
\Delta_1^{g_d}(x) = \frac{x}{2} \, \left( x - g_d(2)\right)
\end{equation*}
and $g_d(2) \geq 2$. This explains the results for $n=1$ at $x=1$.

For $x >0$ we have $\Delta_{1}^{\psi _{d}}\left( x\right) \geq 0 $
if and only if $x \geq 2^d$. Let $n \geq 2$. Table \ref{table:d2} records our results for Tur\'an inequalities for small $d$.
\begin{table}[H]
\begin{tabular}{r|l|l}
\hline
$ d$ & Tur\'an inequality &   proof / verification\\ \hline \hline
$1$ & $\Delta_n^{\func{id}}(x) \geq 0$ for $n>1$ 
and $x \in \mathbb{R}$ & Heim--Neuhauser 2022          \\
$2,3,4$  &  
$\Delta_n^{\psi _{d}}(x) \geq 0$ for $n>1$ and $x \in \mathbb{R}$             
&         H--N 2022   ($n \leq  10^{2}$)  \\
\hline
\end{tabular}
\caption{\label{table:d2}  Tur\'an inequalities in the $d$ aspect        }
\end{table}
In the case $d=5$ a new feature appears (see Figure \ref{d5}).
Let $3 \leq n \leq 100$ then there are exactly two simple positive zeros.
Their
position implies $\Delta_3^{\psi_5}(1)>0$,  $\Delta_4^{\psi_5}(1)<0$
and $\Delta_n^{\psi_5}(1)>0$ for $ 5 \leq n \leq 100$. We expect that this is true
for all $n \geq 5$.
%%%%%%%%%%%%%%%%%%%%%%%%%%%%%%%%%%%%%%%%%%%%%%%%%%%%%%%%%%%%%%%%%
\begin{minipage}{0.45\textwidth}
\includegraphics[width=\textwidth]{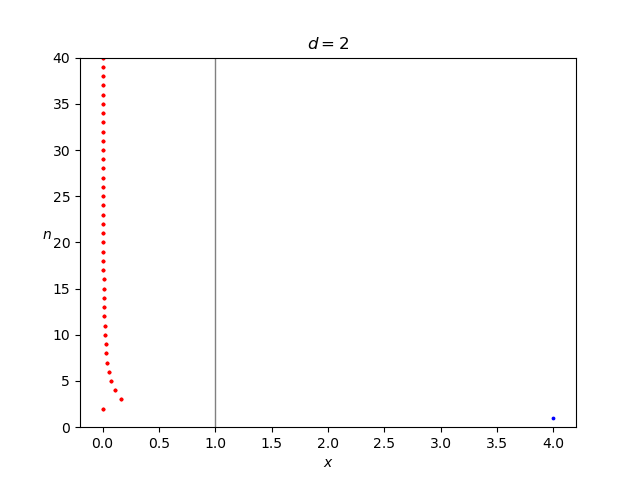}
\end{minipage}
\hfill
\begin{minipage}{0.52\textwidth}
\vspace{1.5cm}
\captionsetup{width=\linewidth}
\captionof{figure}{\label{d2}
The zeros of $\Delta_n^{\psi_2}(x)$ with the largest positive real part for $1 \leq n \leq 40$.
Blue labels the real zeros and red the imaginary zeros.}
\end{minipage}

\begin{minipage}{0.45\textwidth}
\includegraphics[width=\textwidth]{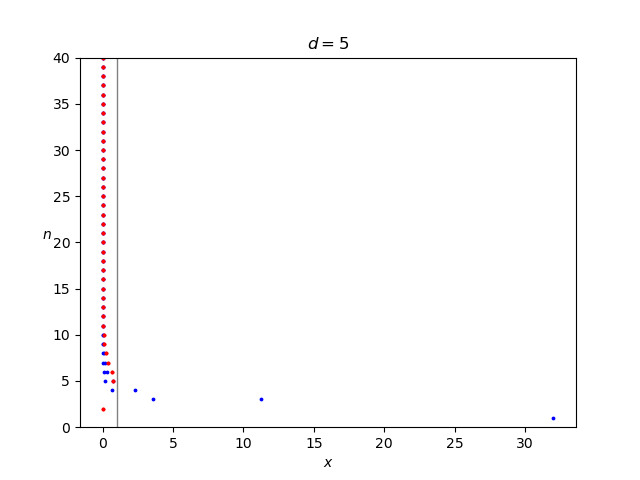}
\end{minipage}
\hfill
\begin{minipage}{0.52\textwidth}
\vspace{1.5cm}
\captionsetup{width=\linewidth}
\captionof{figure}{\label{d5}
The zeros of $\Delta_n^{\psi_5}(x)$ with the largest positive real part for $1 \leq n \leq 40$.
Blue labels the real zeros and red the imaginary zeros.}
\end{minipage}
%%%%%%%%%%%%%%%%%%%%%%%%%%%%%%%%%%%%%%%%%%%%%%%%%%%%%%%%%%%%%%%%%%
%%
%%
%%
%%
\section{Basic
formulas}\label{basic}
Let $g$ be a normalized arithmetic function. Let
$$P_n^g(x):= \frac{x}{n} \sum_{k=1}^n g(k) \, P_{n-k}^g(x), \qquad \text{with }
P_{0}^{g}\left( x\right) =1
.$$
Then $P_n^g(x)$ are polynomials of degree $n$. We refer to \cite{HN20} for a detailed study
of these polynomials. For example
$P_1^g(x)=x$ and $P_2^g(x) = x/2 \, (x +g(2))$. 
\subsection{Coefficients of $P_n^g(x)$}
Let
\begin{equation*}
P_n^g(x) = \sum_{k=0}^n A_{n,k}^g \,\, x^k.
\end{equation*}
Then $A_{0,0}^g=1$. Let $n \geq 1$ then $A_{n,0}^g=0$, $A_{n,1}^g= g(n)\, / \, n$ and $A_{n,n}^g = 1 \, / \, n!$.
We also have \cite{HN20} for $1 \leq m <n$ and $n-m=1,2,3$:
\begin{eqnarray*}
%itemize}
%\item  $
A_{n,n-1}^{g}&=&  \frac{1}{n!} \,\,
g\left( 2\right) \binom{n}{2}
%$
,\\
%\item  $
A_{n,n-2}^{g}&=& \frac{1}{n!} \,\,  \left[
3\left( g\left( 2\right) \right) ^{2}\binom{n}{4}+2g\left( 3\right) \binom{n}{3}\right]
%$
,\\
%\item  $
A_{n,n-3}^{g}&=& \frac{1}{n!} \,\, \left[
15\left( g\left( 2\right) \right) ^{3}\binom{n}{6}+20g\left( 2\right) g\left( 3\right) \binom{n}{5}+6g\left( 4\right) \binom{n}{4}\right]
%$
.
\end{eqnarray*}
%itemize}
%%
%%
\begin{lemma}
Let $g$ be a normalized arithmetic function and 
$$ \Delta_{n}^g(x):= P_n^g(x)^2 - P_{n-1}^g(x) \, P_{n+1}^g(x).$$ 
Then $\Delta_1^g(x) = \frac{x}{2} \left( x-g(2) \right)$ and
$\Delta_2^g(x) =  \frac{x^2}{12} \left( x^2 + 3 \, g(2)^2 - 4 \, g(3)\right)$. Further,
\begin{eqnarray*}
\Delta_3^g(x) & = & \frac{1}{144} x^{6}
 + \frac{1}{48} g\left( 2\right) x^{5}
 + \left( \frac{1}{16}\left( g\left( 2\right) \right) ^{2}
 - \frac{1}{18} g\left( 3\right) \right)  x^{4}
 +
( -\frac{1}{16}\left( g\left( 2\right) \right) ^{3}\\
&&{} + \frac{1}{6} g\left( 3\right) g\left( 2\right)
 - \frac{1}{8} g\left( 4\right)
)  x^{3}
 + \left( -\frac{1}{8} g\left( 4\right) g\left( 2\right)
 + \frac{1}{9}\left( g\left( 3\right) \right) ^{2}\right)  x^{2}.
\end{eqnarray*}
\end{lemma}
This follows from the explicit form of the
polynomials. We have 
\begin{eqnarray*}
P_3^g(x) & = & \frac{x}{6} \left(
x^2 + 3 \, g(2)\, x + 2 \, g(3) \right)
,\\
P_4^g(x) & = & \frac{x}{24} \left(
x^3 + 6 \, g(2) \, x^2 + \left( 8 \, g(3) + 
3 \, g(2)^2 \right) \, x + 6 \, g(4) \right).
\end{eqnarray*}
\subsection{Properties of $\Delta_n^g(x)$}
Let us establish the following notation:
\begin{equation*}
\Delta_n^g(x) = \sum_{k=0}^{2n
} D_{n,k}^g \,\, x^k.
\end{equation*}
In contrast to $P_n^g(x)$,  the coefficients of 
$\Delta_n^g(x)$ are not always non-negative in general.
Nevertheless, we have $\Delta_n^g(0)=0$ and 
the important asymptotic property 
\begin{equation*}
\lim_{x \to \infty} \Delta_n^g(x) = \infty.
\end{equation*}
This follows from $D_{n,2n
} =  \frac{1}{(n!)^2 \, (n+1)}$. Let $n \geq 2$. We can
always factor out $x^2$ and
still have polynomials, since $D_{n,0}^g = D_{n,1}^g =0$. The {\it new constant term}
is given by $D_{n,2}^g$, which does not need to be non-negative:
\begin{eqnarray*}
D_{n,2}^g & = & \frac{1}{n^2} \, \left[
g(n)^2 - \frac{n^2}{n^2-1} \, g(n-1) \, g(n+1) \right] ,
\\
D_{n,3}^g & = & 2\frac{g\left( n\right) }{n^{2}}\sum _{k=1}^{n-1}\frac{g\left( k\right) g\left( n-k\right) }{k}-\frac{g\left( n-1\right) }{n
-1}\sum _{k=1}^{n}\frac{g\left( n+1-k\right) g\left(
k\right) }{2\left( n+1-k\right) k}\\
&&{}-\frac{g\left( n+1\right) }{n
+1}\sum _{k=1}^{n-2}\frac{g\left( n-1-k\right) g\left(
k\right) }{2\left( n-1-k\right) k}. \nonumber
\end{eqnarray*}

\subsection{Special cases}
Let $\{g_d(n)\}$ be a double sequence in $
\mathbb{D}$.
We have $\Delta_1^{g_d}(x) = x \, (x-g_d(2))/2$. Thus, $\Delta_1^{g_d}(x) = 0$ if $x=0$ or $x= g_d(2)$.
Thus, $\Delta_1^{g_d}(x)>0$ if and only if $x \not\in [0, g_d(2)]$. 
This implies that $\Delta_1^{g_d}(x)<0$ for
$x\in \left( 0,g_{d}\left( 2\right) \right) $ and all $d \in \mathbb{N}$.
The case $n=2$ is still directly accessible. We have $\Delta_2^{g_d}(x) = 0$ if $x=0$ or
$x^2= 4\, g_d(3)-3 \,\left( g_{d}\left( 2\right) \right) ^2$. We consider $4\, g_d(3)-3 \,\left( g_{d}\left( 2\right) \right) ^2 \geq 0$ and $x \neq 0$.
Let $g_{d}
= \psi _{d}
$ or $g_{d}=\sigma _{d}$. Then $\Delta_2^{\psi_d}(x) >0$ for $d
\in \mathbb{N}$, especially $\Delta_2^{\psi_1}(x) =x^4 /12$.
%%%%%%%%%%%%%%%%%%%%%%%%%%%%%%%%%%%%%%%%%%%%%%%%%%%%%
%%%%%%%%%%%%%%%%%%%%%%%%%%%%%%%%%%%%%%% Theorem Proof
%%%%%%%%%%%%%%%%%%%%%%%%%%%%%%%%%%%%%%%%%%%%%%%%%%%%%
\section{Proof of Theorem \ref{th:1}}
Our strategy is to utilize the well-known formula (\cite{Ko04}, section 4.7):
\begin{equation*}
P_{n}^{g_{d}}\left( x\right) =
\sum _{k\leq n}\sum _{\substack{m_{1},\ldots ,m_{k}\geq 1 \\ m_{1}+\ldots +m_{k}=n}}
\frac{1}{k!} \, \frac{g _{d}\left( m_{1}\right) 
\cdots g _{d}\left( m_{k}\right) }{m_{1}\cdots m_{k}}  \,\, x^k.
\end{equation*}
\subsection{Lower and upper bounds}
In Section 3 \cite{HN22} we have obtained lower and upper bounds
for $P_n^{\sigma_d}(1)$. The invented proof method can be generalized in a
straightforward manner to obtain the following result for
all double sequences in $
\mathbb{D}$ and the associated polynomials for $x >0 $.
\begin{proposition}\label{prop:1}
Let the double sequence
$\{g_d(n)\}_{d,n \in \mathbb{N}}$ be an element of $
\mathbb{D}$.
Let $n \geq 3$ and $x >0$. 
Then we have for all $d\geq 1$ the following upper and lower bounds. 

Let $n \equiv 0 \pmod{3}$ and $n^{\prime}:= n/3$. Then
\begin{equation*}
\frac{3^{\left( d-1\right) n'}}{\left( n^{\prime}  \right) !} \,\, x^{n^{\prime}}
 \, <    \,       P_n^{g_d}(x)   \,  \leq \, 
3^{\left( d-1\right) n^{\prime}}           P_n^{g_1}(x).
\end{equation*}
Further, let $n \equiv 1 \pmod{3}$ and $n^{\prime}:= (n-4)/3$. Then
\begin{equation*}
\frac{\left( 4 \cdot 3^{n'}\right)^{d-1}}{\left( n^{\prime} \right) !}
\, \left( x^{n^{\prime}+1} + \frac{x^{n^{\prime}+2}}{2} \right) \,
<   \,     P_n^{g_d}(x)         \, \leq \,
\left( 4 \cdot 3^{n^{\prime}}\right)^{d-1}   \,\,    P_n^{g_1}(x).
\end{equation*}
Further, let $n \equiv 2 \pmod{3}$ and $n'
:= (n-2)/3$. Then
\begin{equation*}
\frac{\left( 2\cdot 3^{n'}\right)^{d-1}}{\left( n' \right) !} \,\, x^{n'+1} 
\, <  \, P_n^{g_d}(x) \, \leq \, 
\left( 2\cdot 3^{n'}\right)^{d-1}    P_n^{g_1}(x).
\end{equation*}
Additionally, let $n \equiv 2 \pmod{3}$ and $n \geq 8$. Let $n':= (n-2)/3$. Then
\begin{equation*} \label{improve}
P_n^{g_d}(x) \leq
\frac{\left( 2\cdot 3^{n'}\right) ^{d-1}}{\left( n'\right) !}     \, x^{n'+1}+
\left( 16\cdot 3^{n'-2
}\right) ^{d-1} \, P_n^{g_1}(x).
\end{equation*}
\end{proposition}
%
%
%
%%%%%%%%%%%%%%%%%%%%%%%%%%%%%%%%%%%%%%%%%%%%%%%%%%%%%%%%%%%%%proof verallgemeinert %%%%%%%%%%%%
\subsection{Proof of Theorem \ref{th:1}}
We apply Proposition \ref{prop:1}.
\subsubsection{The case $n\equiv 0 \pmod{3}$}
In the first step we apply Proposition \ref{prop:1}. This leads to
\begin{equation*}
\frac{\left( P_{n}^{g_{d}
}\left( x\right) \right) ^{2}}{P_{n-1}^{g_{d}
}\left( x\right) P_{n+1}^{g_{d}
}\left( x\right) } \geq 
\frac{x^{2n/3}}{\left( \left( n/3\right) !\right) ^{2}P_{n-1}^{g_{1}
}\left( x\right) P_{n+1}^{g_{1}
}\left( x\right)
}\left( \frac{9}{8}\right) ^{d-1}.
\end{equation*}
We choose $r(x)>0$, such that $P_n^{g_1}(x) \leq r(x)^{-n}$ for all $n$. Let 
$$d_{0}=d_{0}\left( n,x\right) =1+\frac{2n}{3\ln \left( 9/8\right) }
\left( \ln \left( n/3\right) -\ln \left( x\right) -3\ln \left(
r\left( {x}\right) \right) \right). $$
Then %the  
$\left\{ P_{n}^{g_{d}}\left( x\right) \right \} _{n}$ is strictly
log-concave at $n$ for $d \geq d_0$, since
\[
\frac{\left( P_{n}^{g_{d}}\left( x\right) \right) ^{2}}{P_{n-1}^{g_{d}}\left( x\right) 
P_{n+1}^{g_{d}}\left( x\right) }\geq \frac{x^{2n/3}}{\left( n/3\right) ^{2n/3}
\left( r\left( {x}\right) \right) ^{-2n}}\left( \frac{9}{8}\right) ^{d-1}>1.
\]
%%%%%%%%%
%%%%%%%%%
%%%%%%%%%
%%%%%%%%%            Rest 1
%%%%%%%%%
%%%%%%%%%
%%%%%%%%%
%%%%%%%%%
\subsubsection{The case $n\equiv 1 \pmod{3}$}
In the first step we apply Proposition \ref{prop:1}. This leads to
\begin{equation*}
\frac{\left( P_{n}^{g_{d}
}\left( x\right) \right) ^{2}}{P_{n-1}^{g_{d}
}\left( x\right) P_{n+1}^{g_{d}
}\left( x\right) } \, \leq \,
\left( \left( \frac{n-1}{3} \right)! \right) ^{2}\,\,
P_n^{g_1}(x)^2 \,\,
x^{-2\left( n-1\right) /3}  \,\, \left( \frac{8}{9}\right) ^{d-1}.
\end{equation*}
We choose $r(x)>0$, such that $P_n^{g_1}(x) \leq r(x)^{-n}$ for all $n$. Let 
%$
\[
d_{0}%^{\prime }
=d_{0}%^{\prime }
\left( n,x\right)
=1+\frac{
2n}{3\ln \left( 9/8\right) }\left( 
\ln \left( \frac{n-1}{3}\right) -3
\ln \left( r\left(
{x}\right) \right) -
\ln \left( x\right)
\right).
\]
Then the sequence
$\left\{ P_{n}^{g_{d}}\left( x\right) \right\} _{n}$ is
strictly log-convex at $n$ for $d \geq d_{0%^{\prime
}$, since
\[
\frac{\left( P_{n}^{g_{d}}\left( x\right) \right) ^{2}}
{P_{n-1}^{g_{d}}\left( x\right) P_{n+1}^{g_{d}}\left( x\right) }\leq \left( \frac{n-1}{3}\right) ^{
2n
/3}\left( r\left(
{x}\right) \right) ^{-2n}x^{2
n
/3}\left( \frac{8
}{9
}\right) ^{d-1}
<1.
\]
%%%%%%%%%%%%%%%%%%%%% Umbruch
\subsubsection{The case $n\equiv 2 \pmod{3}$}
This final case involves some additional considerations. 
Again we first apply Proposition \ref{prop:1} and obtain
\begin{eqnarray*}
&&\left( P_{n}^{g_{d}
}\left( x\right) \right) ^{2}
-P_{n-1}^{g_{d}
}\left( x\right) P_{n+1}^{g_{d}
}\left( x\right)
\\
&\leq &
\left( \frac{\left( 2\cdot 3^{\left( n-2\right) /3}\right) ^{d-1}
x^{\left( n+1\right) /3}}{\left( \left( n-2\right) /3\right) !}+
\left( 16\cdot 3^{\left( n-8\right) /3}\right) ^{d-1}P_{n}^{g_{1}
}\left( x\right) \right) ^{2}
\\
&&
{}-\frac{
\left( 4\cdot 3^{\left( n-5\right) /3}\right) ^{d-1}\left(
x^{\left( n-2\right) /3}+x^{\left( n+1\right) /3}/2\right)
}{
\left( \left( n-5
\right) /3\right) !}\frac{3^{\left( d-1\right) \left( n+1\right) /3}
x^{\left( n+1\right) /3}}{\left( \left( n+1\right) /3\right) !}
\\
&\leq &\left( \frac{\left( 2\cdot 3^{\left( n-2\right) /3
}\right) ^{d-1}x^{\left( n+1\right) /3}}{\left( 
\left( n-2\right) /3\right) !}\right) ^{2}
\left( 1-\frac{\left( n-2\right) /3}{\left( n+1\right) /3}\left( \frac{1}{x}+\frac{1}{2}\right) \right.
\\
&& {}+
2\left( 8
\cdot 3^{
-2
}\right) ^{d-1}x^{-\left( n+1\right) /3}P_{n}^{g_{1}}\left( x\right)
\left( \left( n-2\right) /3\right) !\\
&&{}\left. +\left( \left( \left( n-2\right) /3\right) !\left( 8
\cdot 3^{
-2
}\right) ^{d-1}x^{-\left( n+1\right) /3}P_{n}^{g_{1}}\left( x\right) \right) ^{2}\right).
\end{eqnarray*}
The last inequality can only be not %small
larger than zero if $0<x<2-\frac{12}{n+4}$.
We choose $r(x)>0$ such that $P_n^{g_1}(x) \leq r(x)^{-n}$ for all $n$. 
Then 
\[
\left( \frac{n-2}{3}\right) !\left( \frac{8
}{9
}\right) ^{d-1}x^{-\left( n+1\right) /3}P_{n}^{g_{1}}\left( x\right) \leq \left( \frac{n-2}{3}\right) ^{\left( n-2\right) /3}\left( \frac{8
}{9
}\right) ^{d-1}x^{-\left( n+1\right) /3}
\left( r\left( {x
}\right) \right) ^{-n} <1  \]
for
\[
d>1+\frac{1}{\ln \left( 9/8\right) }\left( \frac{n-2}{3}
\ln \left( \frac{n-2}{3}\right) -\frac{n+1}{3}\ln \left( x\right) -n\ln \left(
r\left( {x}\right) \right) \right) .
\]
Let $d_{0}%^{\prime \prime }%'
=d_{0}%^{\prime \prime }%'
\left( n,x\right) $ be defined by
\begin{eqnarray*}
& &
1+\frac{1}{\ln \left( 9/8\right) }\left( -\ln \left( \frac{n-2}{3n+3}
\left( \frac{1}{x}+\frac{1}{2}\right) -\frac{1}{3}\right) \right. \\
& &
\phantom{1+\frac{1}{\ln \left( 9/8\right) }\left( \right. }\left. {}+\frac{n-2}{3}
\ln \left( \frac{n-2}{3}\right) -\frac{n+1}{3}\ln \left( x\right) -n\ln \left(
r\left( {x}\right) \right) \right) .
\end{eqnarray*}
Then the sequence $\{P_n^{g_d}(x)\}_d$ is
strictly log-convex for all $d \geq d_{0}
%^{\prime \prime }
$.

\section{Tur\'an inequalities}
Let $\{g_d(n)\}$ be a double sequence in $
\mathbb{D}$.
We are interested in finding the
set of positive real numbers, such
that $\Delta_n^{g_d}(x) \geq 0$, with special emphasis on the behavior at $x=1$.

In \cite{HN21} a conjecture for  $\Delta_n^{\sigma_1}(x)$ was %given
stated, which
generalized a conjecture of Chern--Fu--Tang \cite{CFT18} related to integers $x \geq 2$.
The Chern--Fu--Tang conjecture was proven by Bringmann, Kane, Rolen, and Tripp \cite{BKRT21}.
Recently, a second conjecture \cite{HNT21} was proposed for $\Delta_n^{\sigma_2}(x)$.
We have shown for $x=1$, the case of plane partitions, that $\Delta_n^{\sigma_2}(1)>0$ for almost
all $n$. Finally, Ono, Pujahari, and Rolen \cite{OPR22} have proven that $\Delta_n^{\sigma_2}(1)>0$ for all $n \geq 12$.
We now show that $\Delta_{n}^{\psi _{1}}\left( x\right) \geq 0$
for all $x \in \mathbb{R}$. This is the first case
where a full result on Tur\'an inequalities is obtained for a double sequence in $
\mathbb{D}$ with $d$ fixed.

\subsection{$\Delta_n^{\psi_1}(x) \geq 0$}

We have $\psi_1(n)=1$. The polynomials $P_n^{\psi_1}(x)$ had been studied in
\cite{HLN19} and identified with the $\alpha$-associated Laguerre polynomials. We have
$P_n^{\psi_1}(x)= \frac{x}{n}  L_{n-1}^{(1)}(-x)$, where
\begin{equation*}
\sum_{n=0}^{\infty} L_n^{(\alpha)}(x) \, t^n =
 \frac{1}{(1-t)^{
\alpha +1
}}\,\,  \mathrm{e}^{-x \frac{t}{1-t}}, \qquad
\alpha > -1
.
\end{equation*}
The Laguerre polynomials of degree $n$ are given by $L_n(x)= L_n^{(0)}(x)$.
It is known that $\alpha$-associated Laguerre polynomials for $\alpha \geq 0$
satisfy:
\begin{equation*}\label{sum}
\left(L_n^{(\alpha)}(x)\right)^2 - L_{n-1}^{(\alpha)}(x) \,  
L_{n+1}^{(\alpha)}(x) = \sum_{k=0}^{n-1} 
\frac{ \binom{\alpha + n-1}{n-k}}{ n \, \binom{n}{k}} \left(L_k^{(\alpha -1)}(x)\right)^2 >0.
\end{equation*}
These Tur\'an inequalities are not sufficient to prove $\Delta_n^{\psi_1}(x) \geq 0$.
We have to show for all $x \in \mathbb{R}$ that
\begin{equation*}
\Delta_{n}^{\psi_1}\left( x\right) =  \frac{1}{n^2}  
\left( L_{n-1}^{(1)}(-x) \right)^2 -  \frac{1}{n^2-1}  L_{n-2}^{(1)}(-x) 
\,\, L_{n}^{(1)}(-x)\geq 0.
\end{equation*}
Szeg\H{o} \cite{Sz48} proved in 1948, that $L_n^{(\alpha)}(x) / L_n^{(\alpha)}(0)$ satisfies
Tur\'an inequalities, where $ L_n^{(1)}(0)= n+1$. This proves our claim. 
\subsection{Challenges} We propose three open questions.
\subsubsection{Log-concavity of partition and plane partition numbers}
Reprove the results of Nicolas \cite{Ni78} and Ono et al.\ \cite{OPR22} on
the log-concavity of the partition numbers and the plane partition numbers
utilizing the zero distribution of the polynomials $\left\{ P_{n}^{\sigma _{d}}\left( x\right) \right\} $ for $d=1$ and $d=2$.
\subsubsection{Tur\'an inequalities $\Delta_{n}^{\psi_d}(x) \geq 0$}
Based on our numerical investigations on the zeros of $\Delta_n^{\psi_d}(x)$
we believe that it is very likely
for $2 \leq d \leq 4$ that 
$\Delta_n^{\psi_d}(x) \geq 0$ 
for $n \geq 2$ and $x \in \mathbb{R}$. Prove this observation.
%%%%%%%%%%%%%%%%%%%%%%%%%%%%%%%%%%%%%%%%%%%%%%%%%%%%%%%%%%%%%%%%%%%%%%%  n divisible by 3 %%%%%%
\subsubsection{The case $n \equiv 0 \pmod{3}$}
The following problem was presented at the Conference: 
$100$ %y
Years of Mock Theta Functions in 2022 at %the
Vanderbilt University
(organized by Rolen and Wagner). Prove that $\Delta_n^{\psi_d}\left( x\right) \geq 0$
for all $n \equiv 0 \pmod{3}$ and all $d \in \mathbb{N}$.
%%%%%%%%%%%%%%%%%%%%%%%%%%%%%%%%%%%%%%%%%%%%%%%%%%%
%%%%%%%%%%%%%%%%%%%%%%%%%%%%%%%%%%%%%%%%%%%%%%%%%%%
%%%%%%%%%%%%%%%%%%%%%%%%%%%%%%%%%%%%%%%%%%%%%%%%%%%   END
%{\bf Acknowledgments.}

%%%%%%%%%%%%%%%%%%%%%%
%%%%%%%%%%%%%%%%%%%%%%
\end{document}